\documentclass{article}
\usepackage[utf8]{inputenc}
\usepackage{mathrsfs}
\usepackage{amssymb, amsthm, amsmath}
\usepackage{graphicx}
\graphicspath{{Images/}}
\usepackage{caption, subcaption}
\usepackage{authblk}
\newtheorem{Theorem}{Theorem}[section]
\newtheorem{Example}{Example}[section]
\newtheorem{Remark}{Remark}[section]

\title{How a Change to Topology's Union Axiom Affects Continuity}
\date{\today}
\author{Rachel Bergjord}
\author{Matthew Zabka}
\affil{Southwest Minnesota State University}
\begin{document}

\maketitle

\begin{abstract}
The most general definition of a continuous function requires that the preimage of any open set be open.  Thus, to discuss continuity in the abstract, it is necessary to first define a topology, which tells us which sets in a space are open. Such a topology is given by three axioms that describe how the open sets in a topology behave. In this paper we shall consider a change to one of these axioms and determine how this change affects the continuity of a function.
\end{abstract}
% The abstract above should be written for an advanced audience.  Give a short summary of what the paper is about and what you prove.
%----------------------------------------------------------------------------------------
\section{Introduction}\label{sec:intro}
% This first paragraph needs revision.  Talk about why the definition of a topology makes sense for metric spaces.  Then talk about how, in R, only countable unions are necessary.  Then, explain why this motivates our problem.
	The beginning of topology as its own field of mathematics is often considered the 1895 paper \textit{Analysis Situs} by Henri Poincaré~\cite{AnalysisSitus:2009}. In this paper, Poincaré introduced both homology theory and the fundamental group. %What did Poincare talk about in this paper?
	Ten years later, Maurice Fréchet proposed axioms for convergence and in 1910, David Hilbert suggested axioms for neighborhoods of points. Felix Hausdorff, in his 1914 \textit{Grundzüge der Mengenlehre}, gave the axioms for the metric, limit, and neighborhood approaches for general topological spaces~\cite{Grundzuge:1914}. Finally, in 1925, the axioms of topology were defined as they are today by Russian mathematician Pavel Alexandrov~\cite{History}. %More about the history of toplogy and the formulation of its axioms here.

	In this paper, we are primarily concerned with what we refer to as the \textbf{union axiom} of topology, which states \textit{The union of arbitrarily many open sets is open.} This makes intuitive sense in a metric space, for in a metric space, a set $U$ is open if and only if for each $y \in U$, there is a $\delta >0$ such that $B(y,\delta) \subset U$ where $B(y, \delta)$ is an open ball centered at $y$ with radius $\delta$. Then, given a family of open sets $\{U_\iota\}_{\iota \in \Lambda}$, their union is also open in the metric space definition of open. To see this, assume that $y \in \bigcup_{\iota\in\Lambda} U_\iota$, where each $U_\iota$ is open in a metric space $X$.  Then for some $\iota_0$, we have $y\in U_{\iota_0}$. Since $U_{\iota_{0}}$ is open, there is a $\delta >0$ such that $B(y,\delta) \subset U_{\iota_{0}} \subset \bigcup_{\iota\in \Lambda} U_\iota$. Therefore the union of arbitrarily many open sets is open in a metric space. So the union axiom of topology makes sense when we consider a metric space.

	However, contrast this with the usual topology on $\mathbb{R}$, which is given by a basis consisting of all open sets of the form $(a,b)$. In this case, every open set $V$ in this topology can be written as a countable union of open intervals. 
	
    To see this, let $V$ be any open set in $\mathbb{R}$. Then for each $x \in V$, let $I_x$ denote the largest open interval containing $x$ and contained in $V$. Let
	\begin{center}
	    $a_x= \inf \{ a : (a,x) \subset V \}$ \hspace{0.5 cm} and \hspace{0.5 cm} $b_x= \sup \{b : (x,b) \subset V \}$.
	\end{center}
	Then we must have $a_x<x<b_x$.
	Then $I_x = (a_x,b_x)$, so by construction we have $x \in I_x$ and $I_x \subset V$. Then
	\[
	V=\bigcup_{x \in V} I_x
	\]
	Now suppose that two intervals $I_x$ and $I_y$ intersect. Then their union, which is also an open interval, is contained in $V$ and contains $x$. Since $I_x$ is maximal, we must have $(I_x \cup I_y) \subset I_x$, and similarly $(I_x \cup I_y) \subset I_y$. This can happen only if $I_x = I_y$, so any two distinct intervals in the collection $\mathcal{I}=\{I_x\}_{x \in V}$ must be disjoint. We must now show that there are only countably many distinct intervals in the collection $\mathcal{I}$. Since every open interval $I_x$ contains a rational number and since different intervals are disjoint, they must contain distinct rationals. Therefore $\mathcal{I}$ is countable.

	The above also generalizes to $\mathbb{R}^n$. That is, every open set in the usual topology on $\mathbb{R}^n$ can be written as a countable union of basis elements of the usual topology on $\mathbb{R}^n$
	
	As a result, the usual topology on $\mathbb{R}^n$ would be the same if the union axiom stated \textit{The union of \underline{countably} many open sets is open.} This raises the question \textit{Why is it necessary for the axiom to state that the union of \underline{arbitrarily} many open sets be open?}, as well as \textit{What happens in other topologies if we consider this change to the axiom?} In this paper, we shall consider a modified topology with this change to the union axiom. We call a topology in which the union axiom is changed to only require that countable unions must be open a \textbf{countable union topology}. We prove the following:

\begin{Theorem}\label{thm:Michael}
There is a subset $U$ of $\mathbb{R}$ and a topological space $\mathbb{M}$ such that the indicator function $f:\mathbb{M}\to\mathbb{R}$ given by
\[
f(x) = \begin{cases}
0 & \textrm{ if } x\notin U,\\
1 & \textrm{ if } x \in U
\end{cases}
\]
 is continuous if $\mathbb{R}$ is given the usual topology, however $f$ is \textbf{not} continuous if the union axiom of topology is replaced with a weaker axiom that only requires that \textit{countable} unions of open sets be open.
\end{Theorem}

The constructions in this paper are quite simple, however we have not been able to find in the literature any results that consider the changing topology's axioms. This paper is organized as follows: in Section~\ref{sec:background}, we review the definition of a topology and discuss examples, in particular the Michael line. In Section~\ref{sec:union}, we define a countable union topology and prove Theorem~\ref{thm:Michael}.

\section{Background}\label{sec:background}
Topology is the study of continuity. We use a topology's open sets to abstract the closeness of objects by examining open sets. In what follows, we give a very brief review of topology's axioms and continuity, and we refer an interested reader to~\cite{Munkres:2000qr} for further details.

A \textbf{topology} $\mathcal{T}$ on a set $X$ is a collection of subsets of $X$ that satisfy the following axioms:
\begin{enumerate}
    \item $\emptyset$ and $X$ are in $\mathcal{T}$.
    \item The union of arbitrarily many elements of $\mathcal{T}$ is in $\mathcal{T}$.
    \item The intersection of any finite collection of elements of $\mathcal{T}$ is in $\mathcal{T}$.
\end{enumerate}

We call the elements of $\mathcal{T}$ \textbf{open sets} of $X$, and we refer to the second axiom above as the \textbf{union axiom}.

A collection $\mathcal{B}$ of open sets of $\mathcal{T}$ is a \textbf{basis} for the topology $\mathcal{T}$ if and only if $\mathcal{T}$ equals the collection of all possible unions of elements of $\mathcal{B}$. 

\begin{Example}
A basis for the usual topology on $\mathbb{R}$ is $\mathcal{B}_\mathbb{Q}=\{(a,b)|a,b \in \mathbb{Q}\}$.
\end{Example}

\begin{Example}
The \textbf{Michael line} topology on $\mathbb{R}$, denoted by $\mathbb{M}$, is given by the basis $\mathcal{B}_\mathbb{Q} \bigcup \mathbb{R}\setminus\mathbb{Q}$.
\end{Example}

A function $f:X \to Y$ is \textbf{continuous} if and only if for each open subset $V$ of $Y$, the set $f^{-1}(V)$ is an open subset of $X$.

\begin{Example}\label{ex:indicator}
Let $\mathbb{Q}=\{q_1,q_2,q_3,q_4, \dots \}$. Let $s_1,s_2,s_3, \dots$ be a sequence of positive real numbers such that $\sum_{i=1} ^ \infty s_i=a$ for some finite number $a$. Let $U=\bigcup_{i=1} ^ \infty U_i$ where $U_i$ is an open interval about $q_i$ of length $s_i$. The function $f:\mathbb{M} \to \mathbb{R}$ defined by
	\begin{center}  $f(x)=\begin{cases}
    1 & x\in U\\
    0 & x \notin U
    \end{cases}$
    \end{center}
is continuous.

To see this, let $V$ be any open set in the usual topology on $\mathbb{R}$. Then 
\[
f^{-1}(V) = \begin{cases}
U & \textrm{ if } V \textrm{ contains 1 but not 0}\\
\mathbb{R}\setminus U & \textrm{ if } V \textrm{ contains 0 but not 1}\\
\mathbb{R} & \textrm{ if } V \textrm{ contains both 0 and 1}\\
\emptyset & \textrm{ if } V \textrm{ contains neither 0 nor 1}\\
\end{cases}
\]
In all cases, $f^{-1}(V)$ is open in $\mathbb{M}$.
\end{Example}

%\subsection{Michael Line Topology}
%One topology we consider is the Michael line topology, denoted by $\mathbb{M}$. In this topology, open sets are of the form $U \bigcup P$ where $U$ is open in the usual topology on $\mathbb{R}$ and $P$ is any set of irrational numbers. A reader may check that this is a topology.

%---------------------------------------------------------------------------------------
\section{A Countable Union Topology}\label{sec:union}
A \textbf{countable union topology} $\mathcal{T}_C$ on a set $X$ is a collection of open sets of $X$ that satisfy the following axioms:
\begin{enumerate}
    \item $\emptyset$ and $X$ are open.
    \item The union of countably many open sets is open.
    \item The intersection of any finite collection of open sets is open.
\end{enumerate}

We refer to the second axiom above as the \textbf{countable union axiom}.
A collection $\mathcal{B}$ of open sets of $\mathcal{T}_C$ is a \textbf{basis} for the countable union topology $\mathcal{T}_C$ on $X$ if and only if $\mathcal{T}_C$ equals the collection of all possible countable unions of elements of $\mathcal{B} \bigcup \{X\}$. By definition, every basis for a countable union topology yields a countable union topology.

Let $X_C$ be a set with a countable union topology. Let $Y$ be a set with a topology. We say that $f:X_C \to Y$ is \textbf{continuous} if and only if for each open subset $V$ of $Y$, the set $f^{-1}(V)$ is an open subset of $X_C$.

\begin{Remark}
Note that we have defined \underline{continuous} in two different ways:
\begin{enumerate}
    \item to describe a function between topologies, $f:X \to Y$
    \item to describe a function from a countable union topology to a topology, \\ $f:X_C \to Y$.
\end{enumerate}
The version of continuous we use will be clear from the context.
\end{Remark}

\begin{Remark}
In the above definition of continuous, $Y$ may have a either a topology or a countable union topology.
\end{Remark}

\begin{Remark}
Suppose $f:X_C \to Y$ is continuous. If we give $Y$ a countable union topology, we have that the topology on $Y$ is finer than the countable union topology on $Y$. Then every open subset $V$ in the countable union topology on $Y$ is also open in the topology on $Y$. So $f^{(-1)}(V)$ is an open subset of $X_C$ for every $V$ in the countable union topology on $Y$. So $f:X_C \to Y$ where $Y$ has a countable union topology is continuous.
\end{Remark}

\begin{Remark}
Given a basis $\mathcal{B}$, let $\mathcal{T}$ be the topology generated by $\mathcal{B}$, and let $\mathcal{T}_C$ be the countable union topology generated by $\mathcal{B}$. If there is no countable basis $\mathcal{B}'$ that gives the topology $\mathcal{T}_C$, then $\mathcal{T} \neq \mathcal{T}_C$. In particular, $\mathcal{T}_C \subset \mathcal{T}$. Also note that given a topology $Y$, every function that is continuous from $\mathcal{T}_C$ to $Y$ is also continuous from $\mathcal{T}$ to $Y$. However, the converse is not true as we shall show in Section \ref{sec:proof}.
\end{Remark}

%ALSO, given a basis B, what is the result of looking at the topology generated by B and the countable union topology generated by B

% Add a remark here about the relationship between 'continuous' and 'countably continuous'  Maybe also the relationship between $X$ and $X_C$

Recall from Section~\ref{sec:intro} that in the usual topology on $\mathbb{R}$, every open set can be written as a countable union of open sets. Therefore, when we consider the countable union topology on $\mathbb{R}$, it is the same as the usual topology on $\mathbb{R}$. We shall now give a basis whose usual topology is different from its countable union topology.

% A lot more needs to go here.  For one thing, a countable union topology has not yet been defined, so we must define it in general.  In particular, what is the finite union topology generated by a basis?

\subsection{Countable Union Michael Line Topology}
Let $\mathcal{B}_\mathbb{Q}=\{(a,b)|a,b \in \mathbb{Q}\}$ and consider the Michael line's usual topology, given by the basis
    \begin{equation*}
        \mathcal{M}=\mathcal{B}_\mathbb{Q} \cup \{\{x\}| x \in (\mathbb{R} \setminus \mathbb{Q})\}.
    \end{equation*}
Note that $\mathcal{B}_\mathbb{Q}$ is a basis for the usual topology on $\mathbb{R}$. Let $\mathbb{M}$ be the topology generated by $\mathcal{M}$. Then $\mathbb{M}$ is known as the \textit{Michael line topology}. Also note that since $\mathcal{B}_\mathbb{Q}$ is a basis for the usual topology on $\mathbb{R}$, the Michael line topology is finer than the usual topology on $\mathbb{R}$. In particular, every function that is continuous in the usual topology on $\mathbb{R}$ is also continuous in $\mathbb{M}$.

Now consider the countable union topology given by the basis $\mathcal{M}$. Let $\mathbb{M}_C$ denote this topology. Since every irrational singleton is in $\mathcal{M}$, $\mathcal{M}$ is not countable. For this reason, many sets such as $\mathbb{R} \setminus \mathbb{Q}$ are open in $\mathbb{M}$, but not in $\mathbb{M}_C$. This difference also results in functions that are continuous from $\mathbb{M}$ to $\mathbb{R}$, but  are not continuous from $\mathbb{M}_C$ to $\mathbb{R}$.  The leads to the proof of our main result.
\vspace{0.1in}

\subsection{Proof of Theorem~\ref{thm:Michael}}\label{sec:proof}

\begin{proof}
Since the rational numbers are countable, we can write 
\[
\mathbb{Q}=\{q_1,q_2,q_3,q_4, \dots \}.
\]
Let $s_1,s_2,s_3, \dots$ be a sequence of positive real numbers such that $\sum_{i=1} ^ \infty s_i=a$ for some finite number $a$. Let $U_i$ be an open interval about $q_i$ of length $s_i$ and let $U=\bigcup_{i=1} ^ \infty U_i$. Then $U$ is the union of infinitely many open sets, so it is open and $\mathbb{Q} \subset U$. Let $\lambda(A)$ denote the Lebesgue measure of A. Now, we know that $\lambda(U_i)$ is $s_i$, so 
\[
0 < s_i=\lambda(U_i) \leq \lambda(U)=\lambda(\bigcup_{i=1} ^ \infty U_i) \leq \sum_{i=1} ^ \infty \lambda (U_i)=a
\]
In particular, the measure of $U$ is positive and finite. Let $f$ be the indicator function:
\[f=\begin{cases}
    1 & x \in U\\
    0 & x \notin U
\end{cases}
\]
Example~\ref{ex:indicator} shows that the indicator function $f$ is continuous in the Michael line usual topology, $\mathbb{M}$.
Now consider the countable union topology $\mathbb{M}_C$. Then, for $0<\epsilon<1$, $f^{-1}(1-\epsilon, 1+\epsilon)=U$ where $U=\bigcup_{i=1} ^ \infty U_i$. In particular, $U$ is a countable union of open sets so it is open in $\mathbb{M}_C$. Suppose $f^{-1} (-\epsilon,\epsilon)=\mathbb{R} \setminus U$  were open. Then $\mathbb{R} \setminus U$ would be $\bigcup _{i=1}^\infty B_i$ where $B_i$ is either of the form $(a_i,b_i)$ or ${x_i}$ where $x_i$ is irrational. Since $\mathbb{Q} \subset U$, $\mathbb{R} \setminus U$ does not contain any rational numbers, so $B_i$ cannot be of the form $(a_i,b_i)$. Then $\mathbb{R} \setminus U = \bigcup_{i=1} ^\infty x_i$ would be countable. Since $\lambda (U)$ is finite, $\lambda(\mathbb{R} \setminus U)$ is uncountable, so this case is also not possible. Therefore $\mathbb{R} \setminus U$ is an uncountable union of open sets so it is not open. Therefore $f:\mathbb{M} \to \mathbb{R}$ is continuous while $f:\mathbb{M}_C \to \mathbb{R}$ is not.
\end{proof}

Although functions from $\mathbb{R}$ to $\mathbb{R}$ are not affected when we change the domain to a countable union topology, we see a change in this function $f:\mathbb{M} \to \mathbb{R}$ when we give the domain a countable union topology so that $f:\mathbb{M}_C \to \mathbb{R}$. The change causes the set $\mathbb{R} \setminus U$ to no longer be open, which results in our function becoming discontinuous. 

Many discontinuous functions from the usual topology on $\mathbb{R}$ into $\mathbb{R}$ become continuous if the topology of the domain is changed to the Michael line usual topology. When all of the discontinuities of a function in the usual topology on $\mathbb{R}$ occur at irrational points in the domain, the function will be continuous when the domain is given the Michael line usual topology. This is a result of each irrational singleton being open. For this reason, we would expect the function $f$, given above, to be continuous. However, when we only require that countable unions of open sets be open, this is not the case. This function is just one example of a function that we would expect to be continuous in $\mathbb{M}$, but $f$ is not continuous under our countable union axiom, when $f:\mathbb{M}_C \to \mathbb{R}$. This result demonstrates the importance of the union axiom of topology as it is.

\bibliographystyle{unsrt}
\bibliography{bibliography}

\end{document}